\documentclass[reqno,12pt]{amsart}
\usepackage{amscd,amsfonts,mathrsfs,amsthm,amssymb, amsmath}
\usepackage{enumerate}
\usepackage{bbm}
\usepackage{multicol}
\usepackage{color}
\usepackage{array}
\usepackage{ae}
\usepackage{accents}
\usepackage{epsfig}
\usepackage{pstricks, pst-node, pst-plot, pst-math}
\usepackage[e]{esvect}
\usepackage[nobysame]{amsrefs}
\def\y{\textsl{u}}
\def\x{\textsl{x}}
\def\comment#1 {{\color{red}(Comment: #1) }}

\def\real     #1{{\mathbb R^{#1}}}
\def\complex  #1{{\mathbb C^{#1}}}
\def\natural  #1{{\mathbb N^{#1}}}

\def\dt       {\frac{d}{dt}\,}

\newtheorem{theorem}{Theorem}[section]

\newtheorem{lemma}[theorem]{Lemma}

\newtheorem*{thma}{Theorem A}
\newtheorem*{thmb}{Theorem B}
\newtheorem*{thmc}{Theorem C}

\theoremstyle{definition}
\newtheorem{remark}[theorem]{Remark}
\newtheorem{example}[theorem]{Example}

\def\pproof#1{\@ifnextchar[\opargproof
{\opargproof[\it Proof of #1.]}}
\def\opargproof[#1]{\par\noindent {\bf #1 }}

\makeatother

\numberwithin{equation}{section}

\begin{document}

\title[Local non-collapsing]{Local non-collapsing of volume for the Lagrangian mean curvature flow}
\author[K. Smoczyk]{\textsc{Knut Smoczyk}}
\address{%
	Knut Smoczyk\newline
	Institut f\"ur Differentialgeometrie\newline
	and Riemann Center for Geometry and Physics\newline
	Leibniz Universit\"at Hannover\newline
	Welfengarten 1\newline
	30167 Hannover, Germany\newline
	{\sl E-mail address:} {\bf smoczyk@math.uni-hannover.de}
}

\date{}
\subjclass[2010]{Primary 53C44, 53C21, 53C42}
\keywords{%
	Lagrangian mean curvature flow, 
	singularities,
	almost calibrated,
	type-II}
\thanks{Supported by DFG SM 78/6-1}

\begin{abstract} 
	We prove an optimal control on the time-dependent measure of a measurable set under a reparametrized Lagrangian mean curvature flow of almost calibrated submanifolds in a Calabi-Yau manifold. Moreover we give a classification of those Lagrangian translating solitons in $\complex{m}$ that evolve by this reparametrized flow.
\end{abstract}

\maketitle

\section{Introduction}
	Suppose $(N,J,g_N)$ is a Calabi-Yau manifold of real dimension $n=2m$. A smooth immersion $F:M\to N$ is called \textit{Lagrangian}, if $m=\dim M$ and $F^*\omega_N=0$, where $\omega_N$ denotes the K\"ahler form on $N$. Let $dz$ denote the complex volume form on $N$. A Lagrangian immersion is called \textit{special Lagrangian}, if the Lagrangian is calibrated with respect to the real part of the complex volume form, i.e. if $F^*dz$ coincides with the induced volume form $d\mu$ on $M$. 
	
	It is well known, that for general Lagrangian immersions $F:M\to N$ one has
	$$F^*dz=e^{i\phi}d\mu,$$
	with a multi-valued phase function $\phi$. Thus the phase $\phi$ vanishes for special Lagrangians. Since special Lagrangians are calibrated they are volume minimizing in their homology class, in particular the mean curvature vector field $H$ vanishes. In fact, the phase $\phi$ is related to the mean curvature vector field $H$ by $d\phi=\Theta_H$, where $\Theta_H(\cdot)=\langle H,J\cdot\rangle$ denotes the mean curvature $1$-form on $M$. It is well known that $d\Theta_H=0$ and that the cohomology class of $\Theta_H$ coincides up to a multiple of $\pi$ with the first Maslov class $m_1$ of $M$. Any locally defined potential $\alpha$ for $\Theta_H$ is called a Lagrangian angle and must coincide up to a constant with the phase function $\phi$. Conversely, if $M$ is minimal Lagrangian, then $\phi=\phi_0$ for some constant $\phi_0$. Thus minimal Lagrangians are also calibrated (special) with respect to the complex volume form $e^{-i\phi_0}dz$.
	
	A Lagrangian immersion is called \textit{almost calibrated}, if $\cos(\phi-\phi_0)>0$ for some choice of constant phase shift $\phi_0$. Therefore a Lagrangian is almost calibrated, if and only if $\Theta_H=d\alpha$ with some globally defined Lagrangian angle $\alpha$ with $\cos\alpha>0$. Since the Lagrangian angle satisfies the heat equation, it follows immediately from the maximum principle that the condition to be almost calibrated is preserved during the flow.
	
	One of the most important open questions in the (Lagrangian) mean curvature flow is the classification of possible singularities. In general these fall into two categories, type-I and type-II singularities. The first category describes those finite time singularities for which
	$$\limsup_{t\to T}\bigl(|A|^2(T-t)\bigr)$$
	is bounded, where $T$ denotes the singular time and $|A|$ denotes the norm of the second fundamental form. This category is relatively well understood and it follows from the monotonicity formula of Huisken \cite{huisken}  that a rescaled subsequence converges to a self-similarly shrinking solution. In the compact case, one easily observes from the elliptic equation induced for the Lagrangian angle (see \cite{sm habil}), that there do not exist any compact self-shrinkers with trivial Maslov class (so in particular the first Betti number of a compact self-shrinker must be positive). This was later extended by Neves in \cite{neves} to the case of forming type-I singularities (even non-compact) in the zero Maslov class case and he proved that they never occur. Therefore in the Lagrangian mean curvature flow with zero Maslov class, only type-II singularities need to be studied. This applies in particular to the case of almost calibrated Lagrangians, which was first observed by Wang \cite{wang}.
	
	Nevertheless, it is a hard and in general an unsolved question, what these singularities are. From the general theory for type-II singularities it follows that the tangent flow of such singularities will be an eternal Lagrangian mean curvature flow in $\complex{m}$ with uniformly bounded second fundamental form. Possible candidates are translating solitons, special Lagrangians and products of these types. Besides these there might exist various other types of eternal solutions. A very interesting class of translating solitons was found by Joyce, Lee and Tsui in \cite{joyce}, in particular there exist non-trivial translating solitons with arbitrary small oscillation of their Lagrangian angle. Recently Kunikawa \cite{kunikawa} proved that there do not exist non-flat complete Lagrangian eternal solutions with nonnegative Ricci curvature to the almost-calibrated Lagrangian mean curvature flow in $\complex{m}$ with $\cos\phi\ge\epsilon>0$ for some constant $\epsilon$.
	
	To understand how type-II singularties form one needs to take a closer look at the formation itself. Hopefully, one can then exclude certain types of eternal solutions in the corresponding tangent flow. E.g. in \cite{nevestian}, Neves and Tian classified two-dimensional translating solitons under various conditions, one of which is a control on the volume growth. We believe that the following theorem about the \textit{local non-collapsing of volume} under the Lagrangian mean curvature flow will be useful in understanding many aspects related to the volume in more detail, since it gives an optimal control on the time-dependent measure of a measurable set under the Lagrangian mean curvature flow of almost calibrated submanifolds in a Calabi-Yau manifold.
	
	\begin{thma}\label{thm a}
		Let $(N,J,g_N)$ be a Calabi-Yau manifold and suppose $F_0:M\to N$ is a compact almost calibrated Lagrangian immersion and let $\alpha$ denote a choice for the Lagrangian angle with $\cos\alpha>0$. If $F:M\times[0,T)\to N$ is a smooth solution of the reparametrized Lagrangian mean curvature flow
		\begin{equation}\label{eq rlmcf}
		\dt F=H-\tan\alpha\cdot  JH,\quad F(\cdot,0)=F_0,\tag{\text{$\ast$}}
		\end{equation} 
		then for each measurable set $\Omega\subset M$ and any time $t\in[0,T)$ we have
		\begin{equation}
		\frac{1}{\epsilon}\cdot \operatorname{vol}_0(\Omega)\ge \operatorname{vol}_t(\Omega)\ge \epsilon\cdot\operatorname{vol}_0(\Omega),
		\end{equation}
		where $\epsilon$ is the uniform constant given by
		$$\epsilon:=\min_{M\times\{0\}}\cos\alpha>0$$
		and $\operatorname{vol}_t$ denotes the induced measure on $M$ at time $t$.
	\end{thma}
\begin{remark}~\\[-30pt]
	\begin{enumerate}[(i)]
		\item
		The flow described in \eqref{eq rlmcf} differs from the standard Lagrangian mean curvature flow only by a tangential variation, i.e. there exists a time dependent smooth family of diffeomorphisms $\phi_t:M\to M$ such that $\tilde F(x,t):=F(\phi_t(x),t)$ evolves by the usual mean curvature flow. In particular the flow describes the same evolving submanifolds in $N$. 
		\medskip
		\item
		Since the first Maslov class is trivial on an almost calibrated Lagrangian, almost calibrated Lagrangian submanifolds cannot develop singularities of type-I. Hence all possible finite time singularities are of type-II and the tangent flow of such singularities gives eternal solutions of the Lagrangian mean curvature flow in $\complex{m}$ with bounded second fundamental form (the tangent flow of the reparametrized flow in Theorem A will also be the same reparametrized Lagrangian mean curvature flow in $\complex{m}$). Moreover these eternal solutions must be almost calibrated as well (that $\cos\alpha>0$ follows from $\cos\alpha\ge 0$, the real analyticity of the submanifolds and the strong elliptic maximum principle.)
		\item
		As long as $\cos\alpha\ge\epsilon>0$ during the flow, one can drop the compactness assumption in Theorem A and in such situations it equally well holds in the complete case.
	\end{enumerate}
\end{remark}
\begin{example} 
	The \textit{grim reaper} $\Gamma\subset\complex{}$ given by the graph of the function
	$\y:(-\pi/2,\pi/2)\to\real{}$,
	$$\y(\x)=\log\frac{1}{\cos \x}$$
	is a translating Lagrangian soliton, translating with constant speed $1$ in direction of 
	$V:=\partial/\partial\y$. A short computation shows $d\x=\Theta_H$, $V=H+\nabla \y$. So in particular $\alpha:=\x$ is a Lagrangian angle and
	$$\y=\log\frac{1}{\cos\x} $$
	or 
	$$e^\y\cos\x=1.$$
	Since $e^\y\cos\x$ is constant, we get $d\y=\tan\x\cdot  \Theta_H$, so that with $\alpha=\x$ 
	$$V=H+\nabla \y=H-\tan\alpha\cdot JH.$$
	The same holds for the product $\Gamma\times \Sigma$ of $\Gamma$ with a minimal Lagrangian submanifold $\Sigma\subset\complex{m-1}$. Thus these translating solitons evolve by the reparametrized Lagrangian mean curvature flow given in Theorem A.
\end{example}
A special case of the above example is the product of the grim reaper $\Gamma$ with a flat Lagrangian subspace $\Sigma\subset\complex{m-1}$. This translator usually appears as the blow-up model of the type-II singularities forming in the evolution of immersed Lagrangian spheres, e.g. this is the case for a large class of equivariant spheres containing the Whitney spheres (see \cite{savassmoczyk}). Since the reparametrized Lagrangian mean curvature flow naturally appears for some translating solitons and seems to favor them, we want to understand this in more detail. Without loss of generality we may assume that the origin is contained in $M$. The next theorem gives a classification of these translators.
\begin{thmb}
	Let $M\subset\complex{m}$ be a complete translating soliton, $0\in M$, translating in direction of a unit vector $V\in\complex{m}$ and let $\x, \y$ be the two coordinate functions on $M$ induced by the plane $-JV\wedge V$, i.e. $\y(p):=\langle V,p\rangle $ and $\x(p):=\langle -JV,p\rangle$ for any $p\in M$. Then $\alpha:=\x$ is a Lagrangian angle. Moreover, the following statements are equivalent.
	\begin{enumerate}[\rm (a)]
		\item 
		The translator evolves according to \eqref{eq rlmcf}.
		\item
		The function $e^\y\cos\alpha$ admits a local extremum.
		\item
		$e^\y\cos\alpha$ is constant.
		\item
		$M=\Gamma\times\Sigma$, where $\Sigma\subset\complex{m-1}$ is a minimal Lagrangian submanifold and $\Gamma$ is the grim reaper given by the function $\y(\x)=-\log\cos\x$. 	
	\end{enumerate}
\end{thmb}
	All products $M=\Gamma\times\Sigma$ of the grim reaper with a minimal Lagrangian submanifold $\Sigma$ have in common that $\cos\alpha>0$ and $\inf_M\cos\alpha=0$. This implies in particular that these translating solitons \textit{do not occur} as a blow-up of a type-II singularity on a compact almost calibrated Lagrangian since for such compact Lagrangians we have a uniform lower bound $\cos\alpha\ge\epsilon>0$ for all $t\in[0,T)$. On the other hand this argument does not exclude translating solitons of the form $M=\ell\times\Sigma$, where $\ell$ is a straight line in direction of $V$. For those solitons however we observe that the coordinate function $\y$ is unbounded from below in contrast to those given by $\Gamma\times\Sigma$. That this is not a coincidence will be shown in the next theorem.
\begin{thmc}
	Let $M\subset\complex{m}$ be a complete translating soliton, $0\in M$, translating in direction of a unit vector $V\in\complex{m}$ and with bounded second fundamental form. Let the functions $\alpha:=\x$ and $\y$ be defined as in Theorem B and suppose that $\cos\alpha\ge\epsilon>0$ for some constant $\epsilon$. Then the coordinate function $\y$ is unbounded from above and below.
\end{thmc}
Note, that it was shown by Joyce, Lee and Tsui in \cite{joyce} that for any $\epsilon\in(0,1)$ there exist non trivial translating Lagrangian solitons which satisfy $\cos\alpha\ge \epsilon$. 

In Theorem C we impose the boundedness of the second fundamental form to guarantee the Omori-Yau maximum principle is applicable and hence this condition can be relaxed as long as the Omori-Yau maximum principle still holds. On the other hand, the assumption on the boundedness of the second fundamental form is quite natural, because this will be valid for any parabolic blow-up of a type-II singularity of the mean curvature flow.

\section{Basic notations}%
Let $(N,J, g_N)$ be a Calabi-Yau manifold and suppose $F:M\to N$ is a Lagrangian immersion. The differential $dF$ will be considered as a $1$-form with values in the pull-back bundle $F^*TN$, i.e. 
$$dF\in\Omega^1(M,F^*TN).$$
Composing $J$ with $dF$ we obtain another $1$-form 
$$\nu=JdF\in\Omega^1(M,F^*TN).$$
	From the Lagrangian condition we deduce  
	$$\nu\in\Omega^1(M,T^\perp M),$$ 
	where $T^\perp M$ is the normal bundle of $M$ with respect to the immersion. 
	
	The first fundamental form $g$ on $M$ is the metric induced by $F$, i.e.
	$$g=F^*g_N.$$
	By definition, the second fundamental tensor $A$ of $F$ is 
	$$A=\nabla dF,$$
	where we will use $\nabla$ to denote any canonical connection induced by the Levi-Civita connections on $TM$ resp. $TN$. 
	Since it is well known that $A$ is normal, i.e. $A\in\Gamma(T^\perp M\otimes T^*M\otimes T^*M)$, the Lagrangian condition implies that the tri-linear form
	$$h(u,v,w)=\langle A(u,v),Jw\rangle$$
	is fully symmetric. In the sequel, let $e_1,\dots, e_m$ be a local orthonormal frame for the tangent bundle $TM$. From the Lagrangian condition we get that $\nu_1,\dots,\nu_m$ with $\nu_k:=\nu(e_k)$ forms a local orthonormal frame for the normal bundle $T^\perp M$. Taking covariant derivatives of $dF$ resp. $\nu$ one obtains for any $u,w\in TM$ the equations
	\begin{eqnarray}
	(\nabla_{u}dF)(w)&=&\sum_{k=1}^mh(u,w,e_k)\nu_k,\label{struc 1}\\
	(\nabla_u\nu)(w)&=&-\sum_{k=1}^mh(u,w,e_k)dF(e_k).\label{struc 2}
	\end{eqnarray}
\section{Variations of Lagrangian immersions}
	Suppose now that for some $T>0$ we have a time dependent smooth map
	$$F:M\times[0,T)\to N$$
	such that each map 
	$$F_t:M\to N,\qquad F_t(p):=F(p,t)$$
	is a smooth Lagrangian immersion into $N$. 
	Let $\frac{dF}{dt}$ be the velocity vector field along $M$ considered as a section in the pull-back bundle $F^*TN$ over $M$. We can thus define two one-forms $\eta,\tau\in\Omega^1(M)$ by
	$$\eta(w):=\left\langle \frac{dF}{dt},\nu(w)\right\rangle,\qquad\tau(w):=\left\langle \frac{dF}{dt},dF(w)\right\rangle,$$
	where $w\in TM$ is arbitrary. 

	Let us first define the bilinear forms $K, L$ by
	$$K(u,w):=(\nabla_u\eta)(w)+\sum_{k=1}^m\tau(e_k)h(u,w,e_k),$$
	$$L(u,w):=(\nabla_u\tau)(w)-\sum_{k=1}^m\eta(e_k)h(u,w,e_k).$$
	Next we compute the evolution of the 1-form $dF$ under the flow. We get
	\begin{eqnarray}
	\bigl(\nabla_\dt dF\bigr)(w)
	&=&\nabla_w\left(\frac{dF}{dt}\right)\nonumber\\
	&=&\nabla_w\left(\sum_{k=1}^m\eta(e_k)\nu(e_k)+\sum_{k=1}^m\tau(e_k)dF(e_k)\right)\nonumber\\
	&=&\operatorname{trace}\bigl(\nabla_w(\eta\otimes\nu+\tau\otimes dF)\bigr)\nonumber\\
	&=&\sum_{k=1}^m\Bigl((\nabla_w\eta)(e_k)+\sum_{l=1}^mh(w,e_k,e_l)\tau(e_l)\Bigr)\nu_k\nonumber\\
	&&+\sum_{k=1}^m\Bigl((\nabla_w\tau)(e_k)-\sum_{l=1}^mh(w,e_k,e_l)\eta(e_l)\Bigr)dF(e_k)\nonumber\\
	&=&\sum_{k=1}^mK(w,e_k)\nu_k+\sum_{k=1}^mL(w,e_k)dF(e_k),\label{evol dF}
\end{eqnarray}
	From this we can derive the evolution equation for the metric
	\begin{eqnarray}
	\Bigl(\dt g\Bigr)(u,w)
	&=&\Bigl\langle \bigl(\nabla_\dt dF\bigr)(u),dF(w)\Bigr\rangle+\Bigl\langle \bigl(\nabla_\dt dF\bigr)(w),dF(u)\Bigr\rangle\nonumber\\
	&=&\bigl(\nabla_u\tau\bigr)(w)+\bigl(\nabla_w\tau\bigr)(u)-2\sum_{k=1}^m\eta(e_k)h(u,w,e_k)\nonumber\\
	&=&L(u,w)+L(w,u),\label{evol g}
	\end{eqnarray}
	Thus the evolution of the volume form $d\mu$ is given by
	\begin{equation}\label{evol dmu}
	\dt d\mu=\frac{1}{2}\operatorname{trace}_g\Bigl(\dt g\Bigr)d\mu=(d^\dagger\tau-\langle \Theta_H,\eta\rangle)d\mu,
	\end{equation}
	where $d^\dagger\tau$ is defined by $d^\dagger \tau=\operatorname{trace}(\nabla\tau)=\sum_{k=1}^m\bigl(\nabla_{e_k}\tau\bigr)(e_k)$,  $H$ denotes the mean curvature vector field $H=\operatorname{trace}A$ and $\Theta_H$ is the mean curvature $1$-form on $M$ given by $\Theta_H(\cdot)=\langle H,J\cdot\rangle$.	In the next step we compute the evolution equation of the second fundamental form.
	\begin{eqnarray*}
	&&\bigl(\nabla_\dt h\bigr) (u,v,w)\\
	&=&\underbrace{\left\langle\nabla_\dt(\nabla dF)(u,v),\nu(w)\right\rangle}_{=:S}
	+\underbrace{\left\langle(\nabla_u dF)(v),J(\nabla_\dt dF)(w)\right\rangle}_{=:T},
\end{eqnarray*}
	where we have used that $J$ is parallel. For $T$ we compute with \eqref{struc 1} and \eqref{evol dF}
	\begin{eqnarray*}
	T&=&\left\langle(\nabla_u dF)(v),J(\nabla_\dt dF)(w)\right\rangle\\
	&=&\sum_{k=1}^mh(u,v,e_k)L(w,e_k).
\end{eqnarray*}
	To compute $S$ we first need to interchange the covariant derivatives. By definition of the curvature tensor of the pull-back bundle we have
	\begin{eqnarray*}
	S&=&\left\langle\nabla_\dt(\nabla dF)(u,v),\nu(w)\right\rangle\\
	&=&\left\langle\nabla_u\bigl(\nabla_\dt dF\bigr)(v)+R_N\hspace{-2pt}\left(\frac{dF}{dt},dF(u)\right)\hspace{-2pt}dF(v),\nu(w)\right\rangle,
	\end{eqnarray*}
	where $R_N$ denotes the curvature tensor on $N$. Taking into account \eqref{evol dF} and $J\nu(e_k)=-dF(e_k)$, the first term on the RHS simplifies to
	\begin{eqnarray*}
	\left\langle\nabla_u\bigl(\nabla_\dt dF\bigr)(v),\nu(w)\right\rangle=(\nabla_uK)(v,w)+\sum_{k=1}^mh(u,w,e_k)L(v,e_k)
	\end{eqnarray*}
	Combining everything gives
	\begin{eqnarray*}
	\bigl(\nabla_\dt h\bigr) (u,v,w)&=&(\nabla_uK)(v,w)\nonumber\\
	&&+\sum_{k=1}^mh(u,w,e_k)L(v,e_k)+\sum_{k=1}^mh(u,v,e_k)L(w,e_k)\\
	&&+\left\langle R_N\hspace{-2pt}\left(\frac{dF}{dt},dF(u)\right)\hspace{-2pt}dF(v),\nu(w)\right\rangle.
\end{eqnarray*}
The mean curvature form $\Theta_H$ is given by $\Theta_H(u)=\sum_{k=1}^mh(u,e_k,e_k)$. Taking into account \eqref{evol g}, we take the trace in the last evolution equation over $v, w$ and obtain 
	\begin{eqnarray*}
		\nabla_\dt \Theta_H=d(\operatorname{trace}(K))=d\bigl(d^\dagger \eta+\langle\tau,\Theta_H\rangle\bigr),
	\end{eqnarray*}
where we have used that the Lagrangian condition and the K\"ahler identity on $N$ imply that the trace
$$\sum_{k=1}^m\left\langle R_N\hspace{-2pt}\left(\frac{dF}{dt},dF(u)\right)\hspace{-2pt}dF(e_k),\nu(e_k)\right\rangle$$
gives a Ricci curvature and thus vanishes since Calabi-Yau manifolds are Ricci flat.

	From this evolution equation we deduce that the Lagrangian angle $\alpha$, i.e. the potential with $d\alpha=\Theta_H$, evolves according to
	\begin{equation}\label{evol alpha}
	\dt\alpha=d^\dagger\eta+\langle\tau,\Theta_H\rangle.
	\end{equation}
	As above let $dz$ denote the complex volume form on the Calabi-Yau manifold. Since $F^*dz=e^{i\phi}d\mu$ with phase function $\phi$ we must have $\alpha=\phi-\phi_0$ for some constant $\phi_0$. Thus 
	$$F^*(e^{-i\phi_0}dz)=(\cos\alpha +i\sin\alpha)d\mu.$$
	Therefore from \eqref{evol dmu} we get the evolution equation
	\begin{eqnarray*}
	\dt F^*(e^{-i\phi_0}dz)
	&=&\left(i\dt \alpha+d^\dagger\tau-\langle \Theta_H,\eta\rangle\right)F^*(e^{-i\phi_0}dz).
	\end{eqnarray*}
This and \eqref{evol alpha} imply
\begin{eqnarray*}
\dt (\cos\alpha\, d\mu)&=&\Bigl(-\sin\alpha(d^\dagger\eta+\langle\tau,\Theta_H\rangle)+\cos\alpha(d^\dagger\tau-\langle \Theta_H,\eta\rangle)\Bigr)d\mu,\\
\dt (\sin\alpha\, d\mu)&=&\Bigl(\cos\alpha(d^\dagger\eta+\langle\tau,\Theta_H\rangle)+\sin\alpha(d^\dagger\tau-\langle \Theta_H,\eta\rangle)\Bigr)d\mu.
\end{eqnarray*}
If we now choose
$$\eta=\Theta_H,\qquad\tau=\tan\alpha\cdot \Theta_H=-d(\log\cos\alpha),$$
we get
\begin{eqnarray*}
&&-\sin\alpha(d^\dagger\eta+\langle\tau,\Theta_H\rangle)+\cos\alpha(d^\dagger\tau-\langle \Theta_H,\eta\rangle)\\
&=&-\sin\alpha (d^\dagger \Theta_H+\tan\alpha|\Theta_H|^2)\\
&&+\cos\alpha\left(\tan\alpha\cdot d^\dagger \Theta_H+\frac{|\Theta_H|^2}{\cos^2\alpha}-|\Theta_H|^2\right)\\
&=&0.
\end{eqnarray*}
Moreover
\begin{eqnarray*}
	&&\cos\alpha(d^\dagger\eta+\langle\tau,\Theta_H\rangle)+\sin\alpha(d^\dagger\tau-\langle \Theta_H,\eta\rangle)\\
	&=&\cos\alpha (d^\dagger \Theta_H+\tan\alpha|\Theta_H|^2)\\
	&&+\sin\alpha\left(\tan\alpha\cdot d^\dagger \Theta_H+\frac{|\Theta_H|^2}{\cos^2\alpha}-|\Theta_H|^2\right)\\
	&=&\frac{d^\dagger \Theta_H}{\cos\alpha} +\frac{\sin\alpha}{\cos^2\alpha} |\Theta_H|^2=-\Delta(\log{\cos\alpha})
\end{eqnarray*}
We summarize this in the following Lemma.
\begin{lemma}
	Under the reparametrized Lagrangian mean curvature flow
	$$\dt F=H-\tan\alpha\cdot JH$$
	we have
	\begin{eqnarray}
		\dt (\cos\alpha\, d\mu)&=&0,\label{eq evol c}\\
		\dt (\sin\alpha\, d\mu)&=&\Delta\left(\log{\frac{1}{\cos\alpha}}\right)\, d\mu.\label{eq evol s}
	\end{eqnarray}
\end{lemma}
\begin{remark}
	If one first chooses $\eta=\theta_H$, then $\dt(\cos\alpha\, d\mu)=0$, if and only if
	$$d^\dagger(\cos\alpha\cdot\tau-\sin\alpha\cdot \theta_H)=0.$$
	Therefore one observes that $\dt(\cos\alpha\, d\mu)=0$, if and only if 
	$$\cos\alpha\cdot\tau=\sin\alpha\cdot\theta_H+\sigma,$$
	where $d^\dagger\sigma=0$, i.e. where $\sigma$ is a smooth time-dependent co-closed $1$-form.
\end{remark}
From the computations above one observes that the Lagrangian angle evolves under the reparametrized Lagrangian mean curvature flow \eqref{eq rlmcf} according to
$$\dt\alpha=\Delta\alpha+\tan\alpha\,|\nabla \alpha|^2.$$
Equation (\ref{eq evol c}) is the key ingredient to prove Theorem A.

\textbf{Proof of Theorem A.}
The evolution equation for $\cos\alpha$ is
$$\dt\cos\alpha=\Delta\cos\alpha+\cos\alpha\,|\nabla\alpha|^2-\frac{1}{\cos\alpha}|\nabla(\cos\alpha)|^2.$$
Thus the parabolic maximum principle implies
\begin{equation} \label{eq evol e}
\min_{M\times\{t\}}\cos\alpha\ge \min_{M\times\{0\}}\cos\alpha=\epsilon,\quad\text{for all }t\in[0,T).
\end{equation} 
Then for each measurable set $\Omega\subset M$ we get
$$\operatorname{vol}_t(\Omega)=\int\limits_{\Omega\times\{t\}}d\mu\ge\int\limits_{\Omega\times\{t\}}\cos\alpha\, d\mu\overset{\eqref{eq evol c}}{=}\int\limits_{\Omega\times\{0\}}\cos\alpha\, d\mu\ge\epsilon\operatorname{vol}_0(\Omega).$$
Similarly
$$\operatorname{vol}_t(\Omega)
=\int\limits_{\Omega\times\{t\}}d\mu
\overset{\eqref{eq evol e}}{\le}\frac{1}{\epsilon}\int\limits_{\Omega\times\{t\}}\cos\alpha \, d\mu\overset{\eqref{eq evol c}}{=}\frac{1}{\epsilon}\int\limits_{\Omega\times\{0\}}\cos\alpha \, d\mu\le\frac{1}{\epsilon}\operatorname{vol}_0(\Omega).$$
This completes the proof of Theorem A.\hfill{$\square$}
\section{Translating solitons}
In general, the equation for a translating soliton $F:M\to\complex{m}$ for the mean curvature flow is
\begin{equation}\label{eq trans}
H=V^\perp,
\end{equation}
where $H$ denotes the mean curvature vector of $M$, $V$ is a constant vector of unit length in $\complex{m}$ and $V^\perp$ denotes the normal part of $V$ along the submanifold. There exist a number of results for translating solitons, e.g. in \cite{mss}, \cite{nevestian}, \cite{sun1} and \cite{sun2}.

In case of a Lagrangian translating soliton $M\subset\complex{m}$ equation \eqref{eq trans} can be expressed in terms of the mean curvature $1$-form $\Theta_H$,
\begin{equation}\label{eq trans3}
\Theta_H=d\x,
\end{equation}
where $\x$ is the coordinate function
$$\x:=-\langle JV,F\rangle.$$
In particular, translating Lagrangian solitons have trivial first Maslov class and are of gradient type. Moreover,
\begin{equation}\label{eq transv}
\alpha:=\x
\end{equation}
 is a choice for the Lagrangian angle $\alpha$. From this one immediately gets
\begin{equation}\label{eq trans1}
V=H+\nabla \y,
\end{equation}
where $\y$ is the coordinate function
$$\y:=\langle V,F\rangle.$$
$|V|=1$ implies
\begin{equation}\label{eq trans2}
1=|\nabla\alpha|^2+|\nabla \y|^2.
\end{equation}
Next we will compute the Laplacians of various functions. From $\alpha=\x$ we get
\begin{eqnarray}
\nabla^2\alpha&=&-h(\nabla \y,\cdot,\cdot),\label{eq trans4}\\
\nabla^2\y&=&h(\nabla\alpha,\cdot,\cdot).\label{eq trans5}
\end{eqnarray}
Taking traces gives
\begin{eqnarray}
\Delta\alpha+\langle\nabla\alpha,\nabla \y\rangle&=&0,\label{eq trans6}\\
\Delta \y-|\nabla\alpha|^2&=&0.\label{eq trans7}
\end{eqnarray}
With \eqref{eq trans2} we derive from \eqref{eq trans7} that
\begin{eqnarray}
\Delta e^\y=e^\y(\Delta \y+|\nabla \y|^2)=e^\y.
\end{eqnarray}
Let us define the function 
$$f:=e^\y\cos\alpha.$$
Since
\begin{eqnarray}
\nabla f=e^\y(\cos\alpha\nabla \y-\sin\alpha\nabla\alpha)=f\nabla\y-e^\y\sin\alpha\nabla\alpha,
\end{eqnarray}
we obtain
\begin{eqnarray}
\nabla^2f&=&\nabla f\otimes\nabla \y +f(\nabla^2\y-\nabla\alpha\otimes\nabla\alpha)-e^\y\sin\alpha(\nabla^2 \alpha+\nabla\y\otimes\nabla\alpha).\nonumber\\
&&\label{eq trans9}
\end{eqnarray}
Hence taking a trace and using \eqref{eq trans6}, \eqref{eq trans7} we get
\begin{eqnarray}
\Delta f-\langle\nabla f,\nabla \y\rangle=0.\label{eq trans8}
\end{eqnarray}
We want to exploit equation \eqref{eq trans9} even further in case $\cos\alpha>0$. To this end observe that from $1=\frac{1}{\cos^2\alpha}-\tan^2\alpha$ we get
\begin{eqnarray*}
&&\nabla^2\y-\nabla\alpha\otimes\nabla\alpha-\tan\alpha (\nabla^2\alpha+\nabla\y\otimes\nabla\alpha)\\
&=&\frac{1}{\cos^2\alpha}(\nabla^2 \y-\nabla\alpha\otimes\nabla\alpha)\\
&&+\frac{\tan\alpha}{\cos\alpha}\Bigl(\sin\alpha(\nabla\alpha\otimes\nabla\alpha-\nabla^2\y)-
\cos\alpha  (\nabla^2\alpha+\nabla\y\otimes\nabla\alpha)\Bigr).
\end{eqnarray*}
The last line can be substituted using equations \eqref{eq trans4}, \eqref{eq trans5} and this gives
\begin{eqnarray*}
	&&\nabla^2\y-\nabla\alpha\otimes\nabla\alpha-\tan\alpha (\nabla^2\alpha+\nabla\y\otimes\nabla\alpha)\\
	&=&\frac{1}{\cos^2\alpha}(\nabla^2 \y-\nabla\alpha\otimes\nabla\alpha)\\
	&&+\frac{\tan\alpha}{\cos\alpha}\Bigl(h(\cos\alpha\nabla\y-\sin\alpha\nabla\alpha,\cdot,\cdot)+(\sin\alpha\nabla\alpha-\cos\alpha\nabla\y)\otimes\nabla\alpha\Bigr).
\end{eqnarray*}
Therefore we have
\begin{eqnarray*}
	&&\nabla^2\y-\nabla\alpha\otimes\nabla\alpha-\tan\alpha (\nabla^2\alpha+\nabla\y\otimes\nabla\alpha)\\
	&=&\frac{1}{\cos^2\alpha}(\nabla^2 \y-\nabla\alpha\otimes\nabla\alpha)+\frac{\tan\alpha}{f}\Bigl(h(\nabla f,\cdot,\cdot)-\nabla f\otimes\nabla\alpha\Bigr).
\end{eqnarray*}
Combining this with \eqref{eq trans9} implies
\begin{eqnarray}
\nabla^2f
&=&\nabla f\otimes\nabla\y+\frac{f}{\cos^2 \alpha}(\nabla^2\y-\nabla\alpha\otimes\nabla\alpha)\nonumber\\
&&+\tan\alpha\Bigl(h(\nabla f,\cdot,\cdot)-\nabla f\otimes\nabla\alpha\Bigr)\nonumber\\
&=&\frac{1}{f}\nabla f\otimes\nabla f+\tan\alpha\cdot  h(\nabla f,\cdot,\cdot)\nonumber\\
&&+\frac{f}{\cos^2 \alpha}\Bigl(h(\nabla\alpha,\cdot,\cdot)-\nabla\alpha\otimes\nabla\alpha\Bigr).\label{eq trans10}
\end{eqnarray}
\textbf{Proof of Theorem B.}
Let us first mention that since $0\in M$ and $\alpha=\x=-\langle p,JV\rangle$, there exists at least one point $p\in M$ with $\y(p)=0$ and $\alpha(p)=0$. Therefore, if $f$ is constant, this constant clearly is $1$.
	\begin{enumerate}[1.]
		\item 
			We prove (a)$\Leftrightarrow$(c):
		\item[]
			From $V=H+\nabla \y$ we see that $V$ takes the form in (\ref{eq rlmcf}), if and only if 
			$$\nabla \y=-\tan\alpha \cdot JH=\tan\alpha\,\nabla\alpha,$$
			i.e. if and only if
			$$\nabla(e^\y\cos\alpha)=e^\y(\cos\alpha\nabla \y-\sin\alpha\nabla \alpha)=0.$$		
		\item
			The equivalence (b)$\Leftrightarrow$(c) follows from the strong elliptic maximum principle applied to the equation \eqref{eq trans8}.
		\item
			(d)$\Leftrightarrow$(c):
		\item[]
			Since $e^\y\cos\alpha\equiv 1$ on the grim reaper $\Gamma$, this implies that $e^\y\cos\alpha$ must be constant and equal to $1$ as well on the product of $\Gamma$ with a minimal Lagrangian submanifold $\Sigma\subset\complex{m-1}$.  So clearly (d) implies (c). It remains to show that (c) implies (d). 
			
			\noindent
			Let us assume that $e^\y\cos\alpha$ is constant. Since the origin is contained in $M$, this constant is $1$. Thus in particular $\cos\alpha>0$ on $M$. From $\nabla f=0$ with $f=e^\y\cos\alpha$ we first observe
		\begin{equation}\label{eq proof1}
		\nabla \y=\tan\alpha\,\nabla\alpha
		\end{equation}
		and then with \eqref{eq trans2}
		$$\sin^2\alpha|\nabla\alpha|^2=\cos^2\alpha|\nabla \y|^2=\cos^2\alpha(1-|\nabla\alpha|^2)$$
		which implies
		$$|\nabla\alpha|^2=\cos^2 \alpha>0.$$
		In particular the kernel of $\Theta_H=d\alpha$ at each point $p\in M$ is $(m-1)$-dimensional. Let $\mathcal D$ denote the corresponding $(m-1)$-dimensional distribution on $M$ defined by $\mathcal D_p:=\operatorname{ker}\Theta_H|_p$. 
		
		\smallskip\noindent\textbf{Claim:} $\mathcal D$ is parallel. \\
		\textit{Proof.}
		\noindent
		From \eqref{eq trans10} and $\nabla f=0$ we derive
		$$h(\nabla\alpha,\cdot,\cdot)=\nabla\alpha\otimes\nabla\alpha.$$
		With \eqref{eq trans4} and \eqref{eq proof1} we then conclude
		\begin{eqnarray*}
			\nabla^2\alpha
			&=&-h(\nabla\y,\cdot,\cdot)\\
			&=&-\tan\alpha\cdot h(\nabla\alpha,\cdot,\cdot)\\
			&=&-\tan\alpha\cdot\nabla\alpha\otimes\nabla\alpha.
		\end{eqnarray*} 
		Now let $\gamma:[0,1]\to M$ be a smooth curve with $\gamma(0)=p$ and let $W$ be the parallel transport of $W_0\in\mathcal D_p$ along $\gamma$. Then we compute along $\gamma$:
		$$\frac{\partial}{\partial s}\bigl(\Theta_H(W)\bigr)=\bigl(\nabla_{\gamma'}\Theta_H\bigr)(W)=-\tan\alpha\cdot\Theta_H(\gamma')\Theta_H(W).$$
		Since the solution of this ODE with $\Theta_H(W)(0)=0$ is unique, it follows that $\Theta_H(W)(s)=0$ for all $s\in[0,1]$ which implies that $W(\gamma(s))\in \mathcal D_{\gamma(s)}$ for all $s$ and that $\mathcal D$ is invariant under parallel transport.
		This proves the claim.\hfill{$\ast$}
		
		\smallskip\noindent
		Therefore the manifold $M$ splits into a Riemannian product of a curve $\beta$ with an $(m-1)$-dimensional submanifold $\Sigma$. Let $\pi :M\to \Sigma$ denote the natural projection. The tangent space $T_{\pi(p)}\Sigma$ is given by $\mathcal D_p$. Since $J\nabla\alpha=H=V^\perp$, the curve $\beta$ lies in the $(\x,\y)$-plane spanned by $V,JV$. Therefore $\Sigma\subset\complex{m-1}$ is also Lagrangian. If at a point $p\in M$ we choose an ONB $\nu_1.\dots,\nu_m$ of the normal space $T_p^\perp M$ with $\nu_1=H/|H|$, then the trace of $A^{\nu_k}:=\langle A(\cdot,\cdot),\nu_k\rangle$ is given by $-\Theta_H(J\nu_k)$ and hence vanishes for all $k\ge 2$. Therefore from the Lagrangian condition and $\mathcal D_p=\operatorname{ker}\Theta_H|_p$ we see that the mean curvature vector of $\Sigma$ vanishes identically and $\Sigma$ is a minimal Lagrangian submanifold. Since $M=\beta\times\Sigma$, with a minimal Lagrangian submanifold $\Sigma$, we finally conclude that $\beta\subset\complex{}$ must itself be a translating soliton in the plane, wich implies $\beta$ must be the grim reaper $\Gamma$.
	\end{enumerate}
\hfill{$\square$}

The Omori-Yau maximum principle \cite{omori} (later extended by Yau \cite{yau}, see also \cite{prs}) states that if $(M,g)$ is a complete Riemannian manifold with sectional curvatures bounded below, then for every $f\in C^2(M)$ that is bounded above there exists a sequence $(x_k)_{k\in\natural{}}\subset M$  such that
$$f(x_k)\ge \sup_M f-\frac{1}{k},\qquad |\nabla f(x_k)|<\frac{1}{k},\qquad \nabla^2 f(x_k)\le\frac{1}{k} g.$$
If a translator $M\subset\complex{m}$ has bounded second fundamental form, then the Gau\ss\ equations imply that all sectional curvatures are bounded. Hence we may apply the Omori-Yau maximum principle to complete translators with bounded second fundamental form.
Since the function $e^\y$ satisfies the equation
$$\Delta e^\y=e^\y$$
the Omori-Yau maximum principle immediately implies that $\y$ cannot be bounded above on any complete translator with bounded second fundamental form. Theorem C claims that this holds also from below, provided the translator is strictly calibrated in the sense $\cos\alpha\ge\epsilon$ for some constant $\epsilon>0$.

\textbf{Proof of Theorem C.}
	Suppose that $\y$ is bounded below by some constant $c_0$. We will derive a contradiction. Since by assumption $\cos\alpha\ge\epsilon>0$, we can choose a constant $\sigma>1$ such that 
	\begin{equation}\label{eq cos}
	\cos(\sigma\alpha)\ge\frac{\epsilon}{2}.
	\end{equation}
	We set $f_\sigma:=e^\y\cos(\sigma\alpha)$ and compute
	\begin{eqnarray}
	\Delta f_\sigma
	&=&\cos(\sigma\alpha)\Delta e^\y+e^\y\Delta\cos(\sigma\alpha)+2\langle\nabla e^\y,\nabla\cos(\sigma\alpha)\rangle\nonumber\\
	&=&f_\sigma+e^\y(-\sigma\sin(\sigma\alpha)\Delta\alpha-\sigma^2\cos(\sigma\alpha)|\nabla\alpha|^2)\nonumber\\
	&&+2\langle\nabla \y,e^\y\,\nabla \cos(\sigma\alpha)\rangle\nonumber\\
	&\overset{\eqref{eq trans6}}{=}&f_\sigma+\langle\nabla \y,e^\y\,\nabla \cos(\sigma\alpha)\rangle-\sigma^2 f_\sigma|\nabla\alpha|^2\nonumber\\
	&=&\langle\nabla \y,\nabla f_\sigma\rangle-f_\sigma|\nabla \y|^2+f_\sigma-\sigma^2 f_\sigma|\nabla\alpha|^2\nonumber\\
	&=&\langle\nabla \y,\nabla f_\sigma\rangle+(1-\sigma^2)f_\sigma|\nabla\alpha|^2.\label{eq cos2}
	\end{eqnarray}
Now
$$\nabla f_\sigma=f_\sigma\bigl(\nabla \y-\sigma\tan(\sigma\alpha)\nabla\alpha\bigr),$$
so that
$$\sigma^2\tan^2(\sigma\alpha)|\nabla\alpha|^2=\left|\frac{\nabla f_\sigma}{f_\sigma}-\nabla \y\right|^2$$
which in view of $|\nabla \y|^2=1-|\nabla\alpha|^2$ implies
$$\frac{(\sigma^2-1)\sin^2(\sigma\alpha)+1}{\cos^2(\sigma\alpha)}|\nabla\alpha|^2=1+\frac{|\nabla f_\sigma|^2}{f_\sigma^2}-2\left\langle\frac{\nabla f_\sigma}{f_\sigma}, \nabla \y\right\rangle.$$
Because $|\nabla \y|^2\le 1$ and $\sigma^2-1\ge 0$ we can apply the Peter-Paul inequality on the right hand side to obtain
\begin{equation}\label{eq cos3}
\frac{\sigma^2}{\cos^2(\sigma\alpha)} |\nabla\alpha|^2\ge \frac{1}{2}-\frac{|\nabla f_\sigma|^2}{f_\sigma^2}.
\end{equation}
Since by assumption $\y\ge c_0$, we may combine
$$\delta:=\inf_Mf_\sigma\ge \frac{\epsilon }{2}\,e^{c_0}>0$$
with \eqref{eq cos2}, \eqref{eq cos3} to get the estimate
\begin{eqnarray}
\Delta f_\sigma
&\le&|\nabla \y|\cdot|\nabla f_\sigma|+(1-\sigma^2)f_\sigma|\nabla\alpha|^2\nonumber\\
&\le&|\nabla f_\sigma|+\frac{1-\sigma^2}{\sigma^2}\cos^2(\sigma\alpha)\left(\frac{f_\sigma}{2}-\frac{|\nabla f_\sigma|^2}{f_\sigma} \right)\nonumber\\
&\le&|\nabla f_\sigma|+\frac{1-\sigma^2}{\sigma^2}\cos^2(\sigma\alpha)\left(\frac{\delta}{2}-\frac{|\nabla f_\sigma|^2}{\delta} \right).\label{eq cos4}
\end{eqnarray}
Thus there exists a positive constant $C=C(\delta,\sigma,\epsilon)$ such that $|\nabla f_\sigma|\le C$ implies
\begin{equation}\label{eq cos5}
\Delta f_\sigma\le\frac{\delta(1-\sigma^2)}{4\sigma^2}\cos^2(\sigma\alpha)\le\frac{\delta\epsilon^2(1-\sigma^2)}{16\sigma^2}.
\end{equation}
But since $\inf_Mf_\sigma=\delta$, the Omori-Yau maximum principle implies that there exists a sequence $(x_k)_{k\in\mathbb{N}}\in M$ such that
$$f_\sigma(x_k)\le\delta+\frac{1}{k},\qquad|\nabla f_\sigma(x_k)|<\frac{1}{k},\qquad\Delta f_\sigma(x_k)\ge -\frac{1}{k}.$$
In view of \eqref{eq cos5} this is for large enough $k$ a contradiction and hence the function $\y$ cannot be bounded below.
\hfill{$\square$}

\bibliographystyle{alpha}

\begin{bibdiv}
\begin{biblist}
\bib{angenent}{article}{
	author={Angenent, S.},
	title={On the formation of singularities in the curve shortening flow},
	journal={J. Differential Geom.},
	volume={33},
	date={1991},
	number={3},
	pages={601--633},
}
%
\bib{chen}{article}{
	author={Chen, J.},
	author={Li, J.},
	title={Singularity of mean curvature flow of Lagrangian submanifolds},
	journal={Invent. Math.},
	volume={156},
	date={2004},
	number={1},
	pages={25--51},
}
%
\bib{huisken}{article}{
	author={Huisken, G.},
	title={Local and global behaviour of hypersurfaces moving by mean
		curvature},
	conference={
		title={Differential geometry: partial differential equations on
			manifolds},
		address={Los Angeles, CA},
		date={1990},
	},
	book={
		series={Proc. Sympos. Pure Math.},
		volume={54},
		publisher={Amer. Math. Soc., Providence, RI},
	},
	date={1993},
	pages={175--191},
}
%
\bib{joyce}{article}{
	author={Joyce, D.},
	author={Lee, Y.-I.},
	author={Tsui, M.-P.},
	title={Self-similar solutions and translating solitons for Lagrangian
		mean curvature flow},
	journal={J. Differential Geom.},
	volume={84},
	date={2010},
	number={1},
	pages={127--161},
}
%
\bib{kunikawa}{article}{
	author={Kunikawa, K.},
	title={Non existence of eternal solutions to Lagrangian mean curvature flow.},
	journal={arXiv:1611.03594v1},
	date={2016},
}
%
\bib{mss}{article}{
	author={Mart\'\i n, F.},
	author={Savas-Halilaj, A.},
	author={Smoczyk, K.},
	title={On the topology of translating solitons of the mean curvature
		flow},
	journal={Calc. Var. Partial Differential Equations},
	volume={54},
	date={2015},
	number={3},
	pages={2853--2882},
}
%
\bib{neves}{article}{
	author={Neves, A.},
	title={Singularities of Lagrangian mean curvature flow: zero-Maslov class
		case},
	journal={Invent. Math.},
	volume={168},
	date={2007},
	number={3},
	pages={449--484},
}
%
\bib{nevestian}{article}{
	author={Neves, A.},
	author={Tian, G.},
	title={Translating solutions to Lagrangian mean curvature flow},
	journal={Trans. Amer. Math. Soc.},
	volume={365},
	date={2013},
	number={11},
	pages={5655--5680},
}
%
\bib{omori}{article}{
	author={Omori, H.},
	title={Isometric immersions of Riemannian manifolds},
	journal={J. Math. Soc. Japan},
	volume={19},
	date={1967},
	pages={205--214},
}
%
\bib{prs}{article}{
	author={Pigola, S.},
	author={Rigoli, M.},
	author={Setti, A. G.},
	title={Maximum principles on Riemannian manifolds and applications},
	journal={Mem. Amer. Math. Soc.},
	volume={174},
	date={2005},
	number={822},
	pages={x+99},
	issn={0065-9266},
	review={\MR{2116555}},
}
%
\bib{savassmoczyk}{article}{
	author={Savas-Halilaj, A.},
	author={Smoczyk, K.},
	title={Lagrangian mean curvature flow of Whitney spheres},
	journal={Preprint},
	date={2018},
}
%
\bib{sm habil}{book}{
	author={Smoczyk, K.},
	title = {Der Lagrangesche mittlere Kr\"ummungsfluss.},
	pages = {102~p.},
	year = {2000},
	publisher = {Leipzig, Univ. Leipzig (Habil.-Schr.)},
	language = {German},
}
%
\bib{sun1}{article}{
	author={Sun, J.},
	title={Rigidity results on Lagrangian and symplectic translating
		solitons},
	journal={Commun. Math. Stat.},
	volume={3},
	date={2015},
	number={1},
	pages={63--68},
}
%
\bib{sun2}{article}{
	author={Sun, J.},
	title={Mean curvature decay in symplectic and lagrangian translating
		solitons},
	journal={Geom. Dedicata},
	volume={172},
	date={2014},
	pages={207--215},
}
%
\bib{wang}{article}{
	author={Wang, M.-T.},
	title={Mean curvature flow of surfaces in Einstein four-manifolds},
	journal={J. Differential Geom.},
	volume={57},
	date={2001},
	number={2},
	pages={301--338},
}
%
\bib{yau}{article}{
	author={Yau, S.-T.},
	title={Harmonic functions on complete Riemannian manifolds},
	journal={Comm. Pure Appl. Math.},
	volume={28},
	date={1975},
	pages={201--228},
	issn={0010-3640},
	review={\MR{0431040}},
}
%
\end{biblist}
\end{bibdiv}

\end{document}